\documentclass[12pt]{article}

\usepackage{amssymb}
\usepackage{latexsym}
\usepackage{amsfonts}
\newtheorem{thm}{Theorem}[section]

\newtheorem{lem}[thm]{Lemma}
\newtheorem{pro}[thm]{Proposition}
\newtheorem{defn}[thm]{Definition}
\title{Functions on groups and computational complexity}

\author{Jean-Camille Birget\thanks{Research supported in part by NSF 
grant DMS-9970471, and in part by NSERC grant 216872-1999} 
       }
\date{}
\begin{document}
\maketitle

\begin{abstract}
We give some connections between various functions defined on finitely 
presented groups (isoperimetric, isodiametric, Todd-Coxeter radius, 
filling length functions, etc.), 
and we study the relation between those functions and the computational 
complexity of the word problem (deterministic time, nondeterministic time, 
symmetric space). We show that the isoperimetric function can always be
linearly decreased (unless it is the identity map). We present a new proof 
of the Double Exponential Inequality, based on context-free languages.
\end{abstract}

%%%%%%%%%%%%%%%%%%%%%%%%%%%%%%%%%%%%%%%%%%%%%%%%%%%%%%%%
% Section 1
%%%%%%%%%%%%%%%%%%%%%%%%%%%%%%%%%%%%%%%%%%%%%%%%%%%%%%%%
\section{Introduction}

The best studied functions on finitely presented group are the so-called
``filling functions'', in particular the isoperimetric functions and the 
isodiametric functions, and more recently, the filling length.
The significance of the filling functions comes from the following: \\  
(1) The connection between the filling functions of a group presentation 
and the computational complexity of the word problem of that presentation. 
By their very definition, filling functions express the difficulty 
or intricateness of certain aspects of a group presentation; hence they 
are a form of complexity by themselves. 
Moreover, there are strong connections between 
certain filling functions and certain computational complexity functions 
of the word problem. Actually, the computational complexity functions of 
the word problem could be considered as filling functions too. \\   
(2) All the filling functions of a group (including the computational 
complexity functions of the word problem) are algebraic invariants of the 
group, in the sense that if one takes a different {\it finite} presentation 
of the same group the filling functions change only linearly 
(i.e., ``up to big-O''). Some filling functions (e.g., the computational
complexity functions of the word problem) are even more strongly invariant;
they do not depend on the presentation, and are invariants of the group
(``up to big-O'') under change of finite set of generators.
Filling functions can give structural information about a group,  
especially when the filling functions are small. 

\smallskip

The term ``filling function'' is inspired from homotopy transformations, which
``fill the space between'' two objects that can be deformed into each other;
this is Gromov's point of view, which has been very influential \cite{Gromov}. 
The same idea can be represented by different images in different 
contexts. From the point of view of computation, this consists of filling in 
the steps between the input and the output of a computation. 
Generally, any transformation that can be decomposed into (or built up from) 
a set of smaller transformations has filling functions that count various 
steps in the transformation.
Filling functions in this general sense have been known since antiquity
(e.g., the most elementary filling functions are the length, area and volume 
of curves, surfaces, bodies; the number of primes in the prime decomposition 
of an integer is another old example of a filling function).
One of the most interesting new developments is the connection between 
``static'' fillings (like length and area) and ``dynamic'' fillings 
(e.g., space and time complexity of computations).

\bigskip

\noindent {\bf Notation and terminology:} \ 
For an alphabet $A$ we let $A^{-1}$ be a disjoint copy of $A$.  
We write $A^{\pm 1}$ for $A \cup A^{-1}$. We will usually assume that
$|A| \geq 2$.
The {\it free monoid} generated by $A^{\pm 1}$ (i.e., the set of all finite
sequences over the alphabet $A^{\pm 1}$, including the empty sequence $1$) is 
denoted by $(A^{\pm 1})^*$; its elements are called {\it words}.
The {\it length} of a word $w \in (A^{\pm 1})^*$ is denoted by $|w|$.
For a set $X$, we denote the cardinality of $X$ by $|X|$.
For a set of words $R$, we denote the sum of the lengths of the words in $R$
by $\|R\|$ ($= \sum_{r \in R} |r|$).
For a word $w = a_1 \ldots a_n$ with $a_i \in A^{\pm 1}$
for $i = 1, \ldots , n$, we define $w^{-1} = a_n^{-1} \ldots a_1^{-1}$.
For a set $R \subseteq (A^{\pm 1})^*$, we define \ 
$R^{-1} = \{r^{-1} : r \in R\}$, and we define $R^{\pm 1} = R \cup R^{-1}$. 
We denote the {\it free group} over the generating set $A$ by FG$(A)$.
We denote the {\it reduction} in FG$(A)$ by {\sf red}; the word 
obtained from $w \in (A^{\pm 1})^*$ by reduction in FG$(A)$ is {\sf red}$(w)$. 
If $G$ is a group with generating set $A$ and if $x, y \in (A^{\pm 1})^*$ 
are words that represent the same element of
$G$ we write $x =_G y$, and we say that $x$ and $y$ are equivalent in $G$ 
(or equivalent modulo $G$); we call the relation $=_G$ on 
$(A^{\pm 1})^*$ the {\it congruence} defining the group $G$. 
If a group $G$ has a presentation 
$\langle A; R \rangle$ with set of generators $A$ and set of relators 
$R \subset (A^{\pm 1})^*$ then $=_{_{\langle A; R \rangle}}$ denotes the 
the congruence on $(A^{\pm 1})^*$ determined by this presentation.  
The congruence class of a word $w$ modulo $=_G$ 
(i.e., the set $\{ x \in (A^{\pm 1})^* : x =_G w\}$) is denoted by $[w]_G$.
If $x, y \in (A^{\pm 1})^*$ are the same word we write $x = y$ and we say 
that $x$ and $y$ are {\it literally equal} 
(the Russian literature calls this ``graphically equal'').
We carefully distinguish between words ($\in (A^{\pm 1})^*$) and elements of 
the free group FG$(A)$ (elements of a free group are equivalence classes of 
words); if $x, y \in (A^{\pm 1})^*$ are equivalent in the free group we 
write \ $x =_{_{{\rm{FG}}(A)}} y$ \ or \ $x =_{_{\rm{FG}}} y$. 
When used between words, ``='' denotes literal equality.  
By $g^h$ we denote the conjugate ($h^{-1}gh$) of $g \in G$ by $h \in G$.  
If $x, y$ are words $(\in (A^{\pm 1})^*)$ we also use the notation 
$x^y$ for the word $y^{-1}xy$.
For a finite presentations $\langle A; R \rangle$ we always assume that 
the empty word 1 is not a relator $(1 \notin R)$.

%%%%%%%%%%%%%%%%%%%%%%%%%%%%%%%%%%%%%%%%%%%%%%%%%%%%%%%%
% Section 2
%%%%%%%%%%%%%%%%%%%%%%%%%%%%%%%%%%%%%%%%%%%%%%%%%%%%%%%%

\section{Filling functions on groups}
 
\begin{defn} \ 
An {\bf isoperimetric} function of a finite presentation 
$\langle A; R \rangle$ of a group is any function 
$P: \mathbb{N} \to \mathbb{N}$ with the following property:
For every $w \in (A^{\pm 1})^*$ such that $w =_{_{\langle A; R \rangle}} 1$ 
there exists a finite sequence $(r_i : i \in I)$ of relators in $R^{\pm 1}$, 
and a finite sequence $(x_i : i \in I)$ of words in $(A^{\pm 1})^*$ 
such that

\medskip

$w  =_{_{{\rm{FG}}}} \prod_{i \in I} r_i^{x_i}$  \ \ and \ \ 
$|I| \leq P(|w|)$. 
\end{defn}

It is known that this definition is equivalent to stating that every word 
$w \in (A^{\pm 1})^*$ such that \ $w =_{_{\langle A; R \rangle}} 1$ \ has 
a van Kampen diagram with {\it area} $\leq P(|w|)$. This characterization 
motivates the term ``isoperimetric''.

It is also equivalent to saying that every word $w \in (A^{\pm 1})^*$ such
that \ $w =_{_{\langle A; R \rangle}} 1$ \ can be rewritten to $1$ in 
at most $P(|w|)$ ``steps''. To define a 
rewrite step we consider the finite {\it rewrite system} with alphabet 
$A^{\pm 1}$ and set of rules $U_R \cup U_{{\rm FG}}$, where \ 

\smallskip

$U_R = \ \{ u \to v : u,v \in (A^{\pm 1})^*, uv^{-1} \in R^{\pm 1} \}$, 
 
\smallskip

$U_{{\rm FG}} = \ \{ aa^{-1} \to 1, 1 \to aa^{-1} : a \in A^{\pm 1}\}$. 

\smallskip

\noindent But in the isoperimetric function we do not count applications 
of the rules in $U_{{\rm FG}}$ (those are automatically part of any group
and are taken for granted; moreover, they don't show up in the van Kampen
diagrams). Note that this rewrite system is {\it symmetric}, 
i.e., if $x \to y$ is a rule then $y \to x$ is  rule too.
The rewriting characterization is the oldest explicit definition of the 
isoperimetric function, under the name ``rewrite distance'' 
(Madlener and Otto \cite{MO}). Gromov \cite{GromovHyper} introduced the 
geometric point of view on these functions; see also Gersten's 
explanations \cite{GerstenIso}.
In \cite{MO} it is proved that the rewrite distances (i.e., 
the isoperimetric functions) are invariants of the group, up to big-O, in
the following sense:  If one takes a different {\it finite} presentation 
of the same group, the isoperimetric function $P(n)$ becomes $c_1 P(c_2n)$, 
where $c_1, c_2 > 0$ are ``constants'' (i.e., they do not depend on $n$, but 
they depend on the two presentations). This fact was rediscovered a little 
later independently by several authors.

Yet another equivalent definition of the isoperimetric function can be 
obtained by using the Cayley 2-complex of the presentation. 
An isoperimetric function of the presentation is any function 
$P: \mathbb{N} \to \mathbb{N}$ with the following property:
For every $w \in (A^{\pm 1})^*$ such that $w =_{_{\langle A; R \rangle}} 1$ 
there exists a 
combinatorial homotopy in the Cayley 2-complex, starting with a loop labeled
by $w$ (with base point $1$) and ending with the trivial loop consisting of
the base point $1$; the area covered by the homotopy 
transformation is $\leq P(|w|)$ (the ``area covered'' consists of all the 
faces used, {\it with multiplicities}, i.e., a same face that is counted 
repeatedly if it is used repeatedly).

\begin{defn} \  
A {\bf filling length} function of a finite presentation 
$\langle A; R \rangle$ of a group is any function
$F: \mathbb{N} \to \mathbb{N}$ with the following property: Every 
word $w \in (A^{\pm 1})^*$ such that $w =_{_{\langle A; R \rangle}} 1$ 
can be rewritten to $1$ (i.e., 
$w \to w_{N-1} \to \ldots \to w_i \to \ldots w_1 =1$),
using the rewrite system
$(A^{\pm 1}, U_R \cup U_{{\rm FG}})$ defined above,
in such a way that the intermediate words in
the rewrite sequence are all of length \ $|w_i| \leq F(|w|)$.
\end{defn}
Another characterization of the filling length is that every word
$w \in (A^{\pm 1})^*$ such that $w =_{_{\langle A; R \rangle}} 1$ 
has a van Kampen diagram which can be contracted to one point, 
homotopically, in such a way that all
intermediary van Kampen diagrams encountered during this contraction have
perimeter $\leq F(|w|)$; see \cite{GerstenRiley1}.

The concept of filling length was introduced by Gromov \cite{Gromov},
and extensively studied by Gersten and Riley \cite{GerstenRiley1}.
They proved that the filling length changes only up to ``big-O''
(just like the isoperimetric function) when the finite presentation is 
changed. 
 
\begin{defn}
\label{isodiametric} \  
An {\bf isodiametric} function of a finite presentation 
$\langle A; R \rangle$ of a group $G$ is any function 
$D: \mathbb{N} \to \mathbb{N}$ with the following property:
For every $w \in (A^{\pm 1})^*$ such that 
$w =_{_{\langle A; R \rangle}} 1$ there exists a finite
sequence $(r_i : i \in I)$ of relators in $R^{\pm 1}$, and a finite sequence
$(x_i : i \in I)$ of words in $(A^{\pm 1})^*$ such that 

\medskip

$w  =_{_{{\rm{FG}}}} \prod_{i \in I} r_i^{x_i}$ \ \ and \ \    
$|x_i| \leq D(|w|)$ \ (for all $i \in I$).
\end{defn}

This function was explicitly introduced by Gersten (see \cite{GerstenIso} 
for references). The above definition is equivalent to saying 
that every word $w \in (A^{\pm 1})^*$ such that 
$w =_{_{\langle A; R \rangle}} 1$ has a van Kampen 
diagram with diameter $\leq D(|w|)$ (the diameter being measured from the 
origin of $w$ on the perimeter of the diagram). 

\medskip

Another characterization of the isodiametric functions is by means of 
a loop complex. For any positive integer $j$ we define a labeled directed 
graph $\Lambda_j$ with labels in $A$, as follows: \\ 
- First, we create a new vertex $v_0$ (the ``origin''). \\  
- Second, for every relator $r \in R$ and every reduced word 
$u \in (A^{\pm 1})^*$ such that $|u| \leq j$, we create a path labeled 
by $u$, starting from vertex $v_0$; at the non-$v_0$ end of this path we 
attach a loop labeled by $r$; in doing this, we create $|u| + |r| -1$ new 
vertices (for every pair $(u,r)$). 
All the paths (and all the loops) are disjoint, except that all the paths 
have the common vertex $v_0$. \\  
- In order to get by with only the alphabet $A$ we replace each edge 
$v_1 \stackrel{a^{-1}}{\to} v_2$ by $v_1 \stackrel{a}{\leftarrow} v_2$,
where $v_1, v_2$ are vertices and $a \in A$. \\  
- We can turn the graph $\Lambda_j$ into a 2-complex by adding a face for 
each loop (the boundary of each face being the corresponding loop). We call 
this the {\bf loop-complex} of radius $j$ of the presentation. \\  
- We can turn the labeled graph $\Lambda_j$ into a 
{\it nondeterministic finite automaton} (an ``{\bf NFA}'') over the alphabet 
$A^{\pm 1}$.  We take the vertex $v_0$ as both start and accept state, 
and we ``symmetrize'' each edge: for every edge 
$v_1 \stackrel{a}{\to} v_2$ we also introduce the ``inverse edge'' 
$v_2 \stackrel{a^{-1}}{\to} v_1$ into the NFA; no new vertices
are created in this symmetrization. 
See \cite{HU} for definitions and basic facts about NFAs.
The language accepted by an NFA $N$ is denoted by $L(N)$.

We call this NFA ``$\Lambda_j$'' too. The context will always make it clear 
whether we refer to the graph, the 2-complex, or the NFA; the vertex set is
the same in the three cases. 

\smallskip

The number of vertices of $\Lambda_j$ is \  
$\leq (\|R\| + |R| \cdot j)  \, (2|A|)^j \leq c^j$, 
where $c > 1$ is a
constant that depends on the presentation $\langle A; R \rangle$.

\medskip

Characterization of the isodiametric functions:
A function $D(\cdot)$ is an isodiametric function of $\langle A; R \rangle$ 
iff for every word 
$w \in (A^{\pm 1})^*$ of length $\leq n$ we have:  \ 
$w =_{_{\langle A; R \rangle}} 1$ iff there is a 
loop in $\Lambda_{D(n)}$, starting and ending at $v_0$, and labeled by a word 
equivalent (in FG$(A)$) to $w$.
This characterization of the isodiametric functions can be reformulated 
as follows: 
\begin{lem}
\label{isodiamNFA} \  
A function $D(\cdot)$ is an isodiametric function of $\langle A; R \rangle$ 
iff for every word $w \in (A^{\pm 1})^*$ of length $|w| \leq n$ we have:

\smallskip

 \ \ \ \ \  \ \ \ \ \  \ \ \ \ \  
 $w =_{_{\langle A; R \rangle}} 1$ \ \ \ iff \ \ \   
 $[w]_{{\rm{FG}}} \cap L(\Lambda_{D(n)}) \neq \emptyset$. 
\end{lem}

The Lemma above can be refined to obtain an upper bound on the minimum 
isoperimetric function $P_{\rm min}(\cdot)$ of $\langle A; R \rangle$.
Lemma \ref{isoperimNFA} is important because it connects $P_{\rm min}$ 
and $D_{\rm min}$; from this we will be able to re-prove the famous
``double exponential inequality'' in a later section.

\begin{lem}
\label{isoperimNFA} \  
Let $\langle A; R \rangle$ be a finite presentation with minimum isoperimetric
function $P_{{\rm min}}$ and isodiametric function $D$. Suppose $\ell(n)$ 
has the property that for all words $w \in (A^{\pm 1})^*$ of length $\leq n$ 
such that $w =_{_{\langle A; R \rangle}} 1$, we have:

\medskip
 
 \ \ \ \ \  \ \ \ \ \  \ \ \ \ \   
 $[w]_{{\rm{FG}}} \cap L(\Lambda_{D(n)})$ \ contains a word of length \   
 $\leq \ell(n)$. 

\medskip

\noindent Then \ $P_{\rm min}(n) \leq \ell(n)$. 
\end{lem}
{\bf Proof:} \ Since $w =_{_{\langle A; R \rangle}} 1$ and $|w| \leq 1$, 
Lemma \ref{isodiamNFA} implies that  \  
$[w]_{{\rm{FG}}} \cap L(\Lambda_{D(n)}) \neq \emptyset$.
Let $z \in [w]_{{\rm{FG}}} \cap L(\Lambda_{D(n)})$ with 
$|z| \leq \ell(n)$. Then $z =_{{\rm{FG}}} w$ and 
$z = \prod_{i \in I} r_i^{x_i}$ (literal equality) for some sequence of 
relators $r_i \in R^{\pm 1}$ and some sequence of words $x_i$ with 
$|x_i| \leq D(n)$. Since the equality \ $z = \prod_{i \in I} r_i^{x_i}$ \  
is literal, and since all $r_i$ are non-empty, we conclude that 
$|I| \leq |z| \leq \ell(n)$. Hence, by the definition of the isoperimetric 
function, $P_{{\rm min}}(n) \leq |I|$.   \ \ \ $\Box$

\bigskip

For later use we introduce a slightly more compact NFA which can play the 
same role as $\Lambda_{D(n)}$ in Lemmas \ref{isodiamNFA} and 
\ref{isoperimNFA} (although it does not accept the same language).
The {\bf tree NFA} of radius $j$ for a presentation $\langle A; R \rangle$ 
is defined as follows: \\ 
- First, we take the Cayley graph of FG$(A)$, truncated to radius $j$ around
the origin; this is a tree of depth $j$, with \ 
$1 + |A^{\pm 1}| \, (|A^{\pm 1}| - 1)^{j-1}$ \ vertices. \\ 
- Second, for every $r \in R$, at every vertex of the above graph we attach 
a loop labeled by $r$. \\ 
- We pick the root vertex of the tree as start and accept state. \\ 
- We symmetrize all edges, in the same way as we did for the NFA 
   $\Lambda_j$. \\  
We call the resulting NFA ``{\it tree}$\Lambda_j$''. 

When $R \neq \emptyset$, the number of states of {\it tree}$\Lambda_j$ is \ 
$\leq \|R\| \cdot (1 + |A^{\pm 1}| \, \sum_{i=0}^{j-1}(|A^{\pm 1}| -1)^i)$. 
Thus, when $R \neq \emptyset$,
the number of states is \ $ \leq \ \|R\| \, (2 \, |A|)^j$.

\begin{lem}
\label{isoperimTreeNFA} \  
Let $\langle A; R \rangle$ be a finite presentation.

\noindent $\bullet$ A function $D(\cdot)$ is an isodiametric function 
of $\langle A; R \rangle$
iff for any word $w \in (A^{\pm 1})^*$ of length $|w| \leq n$ we have:

\smallskip

 \ \ \ \ \  \ \ \ \ \  \ \ \ \ \    
$w =_{_{\langle A; R \rangle}} 1$ \
iff \ $[w]_{{\rm{FG}}} \cap L({\it tree}\Lambda_{D(n)}) \neq \emptyset$.

\medskip

\noindent $\bullet$ Let the minimum isoperimetric function of the presentation 
be $P_{{\rm min}}$ and let  $D$ be an isodiametric function. Suppose $\ell(n)$ 
has the property that for all words $w \in (A^{\pm 1})^*$ of length 
$\leq n$ such that $w =_{_{\langle A; R \rangle}} 1$, we have:

\medskip

 \ \ \ \ \  \ \ \ \ \  \ \ \ \ \    
$[w]_{{\rm{FG}}} \cap L({\it tree}\Lambda_{D(n)})$ \ contains a 
word of length \ $\leq \ell(n)$.

\medskip

\noindent Then $P_{\rm min}(n) \leq \ell(n)$.
\end{lem}
{\bf Proof.} \ The proof is very similar to the proofs of Lemmas 
\ref{isodiamNFA} and \ref{isoperimNFA}. \ \ \ $\Box$

\bigskip

The loop complex  $\Lambda_{D(n)}$ can be {\bf folded}: 
Suppose in the graph we have \
$v_1 \stackrel{a}{\leftarrow} v_2 \stackrel{a}{\to} v_3$,
(or we have \ $v_1 \stackrel{a}{\to}  v_2 \stackrel{a}{\leftarrow} v_3$), 
where $a \in A$ and where $v_1, v_2, v_3$ are vertices. Then we ``glue'' 
$v_1$ and $v_3$ together; 
i.e., we replace $v_1, v_3$ by a new vertex $v_{1,3}$, and we replace 
the two edges $v_1 \stackrel{a}{\leftarrow} v_2$  and 
$v_2 \stackrel{a}{\to} v_3$ by the single edge
$v_{1,3} \stackrel{a}{\leftarrow} v_2$
(respectively, replace $v_1 \stackrel{a}{\to}  v_2$ and 
$v_2 \stackrel{a}{\leftarrow} v_3$ by $v_{1,3} \stackrel{a}{\to} v_2$).
All in-edges (or out-edges) of $v_1$ and $v_3$
become in-edges (respectively out-edges) of $v_{1,3}$.  
This process continues as long as possible. Finally, to obtain a 2-complex
we attach faces on all loops labeled by relators.  

Based on the folded graph one can obtain a ``folded NFA'', called 
f$\Lambda_{D(n)}$, by symmetrizing all edges, and taking $v_0$ as start 
and accept state. The folded NFA is actually deterministic (it is a 
``DFA''), and it has the following property:
 
\begin{lem}
\label{foldedDFA} \  
Let $w \in (A^{\pm 1})^*$ be any word of length $\leq n \in (A^{\pm 1})^*$. 
Then we have $w =_{_{\langle A; R \rangle}} 1$ iff  \ the reduced word 
{\sf red}$(w)$ is accepted by the folded DFA {\rm f}$\Lambda_{D(n)}$.
\end{lem}
{\bf Proof:} \ Let $(N_i : i = 0, 1, \ldots, m)$ be the sequence of NFAs 
obtained by successively folding edges; $N_0 = \Lambda_{D(n)}$, and 
$N_{i+1}$ is obtained from $N_i$ by folding one pair of edges, for 
$i = 0, 1, \ldots, m$; the number of folding steps $m$ is less than the 
number of edges of $\Lambda_{D(n)}$, hence $m \leq c^{D(n)}$ for some 
constant $c > 1$.
 
We want to show that for each NFA $N_i$ ($i = 0, 1, \ldots, m)$, 
the language accepted satisfies \ 
{\sf red}$(L(N_i)) =$ {\sf red}$(L(\Lambda_{D(n)}))$. 
In other words, although $L(N_i)$ 
changes, the set of reductions of all words in $L(N_i)$ does not change.

\medskip

\noindent (1) \ $L(\Lambda_{D(n)}) \subseteq L(N_i)$, \ hence we also 
  have \ 
  {\sf red}$(L(\Lambda_{D(n)})) \subseteq$ {\sf red}$(L(N_i))$ \ 
  for all $i$. 

\smallskip

\noindent This is straightforward, since folding does not destroy any 
reachabilities, but adds additional reachabilities. So, $N_i$ contains 
all the accepting paths of $\Lambda_{D(n)}$ (up to changes of vertex 
names).

\medskip

\noindent (2) \ {\sf red}$(L(N_i)) \subseteq$ 
                {\sf red}$(L(\Lambda_{D(n)}))$ \ for 
all $i$.

\smallskip

\noindent To prove this we use induction on $i$. 
Inclusion (2) is obvious when $i = 0$. 
Suppose {\sf red}$(L(N_{i-1})) \subseteq$ {\sf red}$(L(\Lambda_{D(n)}))$. 
Let $p$ be an accepting path in $N_i$. 
If $p$ does not use an edge involved in the folding step that leads
from $N_{i-1}$ to $N_i$, the path $p$ occurs in $N_{i-1}$ too; so, any 
word accepted by $N_i$ by means of $p$ is also accepted by $N_{i-1}$. 
If $p$ uses one edge involved in the folding step, the names of one of 
the vertices in $p$ changed when $N_{i-1}$ was transformed to $N_i$, but 
the edge labels in $p$ do not change; so, here too, any word accepted 
by $N_i$ by means of $p$ is also accepted by $N_{i-1}$. 
Finally, if two edges in $p$ are folded together then $p$ has the form 
$p_1p_2$, where 

\smallskip

$p_1 \ (v_1 \stackrel{a}{\leftarrow} v_2) \, $
  $(v_2 \stackrel{a}{\to} v_3) \ p_2$

\smallskip

\noindent is a path in $N_{i-1}$, $a \in A$ (we only consider one of
the folding cases; the other is very similar). 
Let $x_1, x_2 \in (A^{\pm 1})^*$ be the labels of $p_1$, respectively 
$p_2$. Along the path $p$, $N_i$ accepts $x_1x_2$, whereas $N_{i-1}$ 
accepts $x_1 aa^{-1} x_2$; but 
{\sf red}$(x_1x_2) =$ {\sf red}$(x_1 aa^{-1} x_2)$.
 \ \ \ $\Box$

\bigskip

We will now define the ``folded'' versions of the above functions.

\begin{defn} \  
A {\bf folded isoperimetric} function of a finite presentation
$\langle A; R \rangle$ of a group is any function
$p: \mathbb{N} \to \mathbb{N}$ with the following property:
Every $w \in (A^{\pm 1})^*$ such that $w =_{_{\langle A; R \rangle}} 1$ 
has a van Kampen diagram 
whose folded area is $\leq p(|w|)$.  The folded area of a van Kampen 
diagram is the number of faces in the 2-complex obtained by folding
the van Kampen diagram.
(Note:  Faces that have the same boundary loop in the 2-complex are
viewed as the same face.)  
\end{defn}

\begin{defn} \ 
A {\bf folded filling length} function of a finite presentation
$\langle A; R \rangle$ of a group is any function
$f: \mathbb{N} \to \mathbb{N}$ with the following property:
Every $w \in (A^{\pm 1})^*$ such that $w =_{_{\langle A; R \rangle}} 1$ 
has a {\em folded} van 
Kampen diagram that admits a homotopy transformation which starts with the 
loop $w$, and ends with the origin point, and with all intermediate loops
of length $\leq f(|w|)$.
\end{defn}
As usual, ``length'' means length of a curve (or a path); it is not just 
the number of different edges; repetitions of edges are counted too.

\begin{defn} \ 
A {\bf folded isodiametric} function of a finite presentation
$\langle A; R \rangle$ of a group is any function
$d: \mathbb{N} \to \mathbb{N}$ with the following property:
Every $w \in (A^{\pm 1})^*$ such that $w =_{_{\langle A; R \rangle}} 1$ 
has a {\em folded} van 
Kampen diagram of diameter (measured from the origin) $\leq d(|w|)$. 
\end{defn}

The folded isoperimetric function and the folded filling length function
seem not to have appeared in the literature. The folded isodiametric 
function has been used a number of times; we will prove later that the 
minimum folded isodiametric function is equal to the minimum isodiametric 
function (and similarly for the filling length function). 

\bigskip

Another function on finite presentations of groups can be defined by 
using the radius of the partial Cayley graphs; these partial Cayley 
graphs are constructed by the following version of the Todd-Coxeter 
process, used here for the word problem for words of length $\leq n$.
We closely follow \cite{Epst} (p.\ 110); see also \cite{StalWolf} (which 
presents a somewhat different graphical version of Todd-Coxeter, however). 

In the process below we use the following definition. In a graph
with origin $v_0$, a {\it hair} is an edge $e$ such that one end vertex 
$v$ of $e$ has the following properties: $v$ has degree 1, and 
$v \neq v_0$.

\bigskip

\noindent {\it Process TC on input} $\langle A; R \rangle$:  

\medskip

create a vertex $v_0$ (called ``origin'');   

{\tt repeat}   

\makebox[.3in][l]{{\footnotesize 1.}} 
             {\tt for} every vertex $v$ of the graph constructed so far: 

\makebox[.6in][l]{} {\tt for} every letter $a \in A$ such that $v$ does 
                                not have an out-edge with label $a$:  
 
\makebox[1.in][l]{}           create a new vertex $(v,a)$ and a new edge 
                                         $v \stackrel{a}{\to} (v,a)$;    
 
\makebox[.6in][l]{} {\tt for} every letter $a \in A$ such that $v$ does 
                                not have an in-edge with label $a$:  
 
\makebox[1.in]{ }    create a new vertex $(v,a^{-1})$ and a new edge 
                                 $v \stackrel{a}{\leftarrow} (v,a^{-1})$;   
 
\makebox[.3in][l]{{\footnotesize 2.}}  
          {\tt for} every vertex $v$ of the graph constructed so far: 

\makebox[.6in][l]{} {\tt for} every relator $r \in R$ which 
                           does not label a loop originating at $v$:  

\makebox[1.in]{}    create a new loop labeled by $r$ and 
                                originating at $v$;   

\makebox[.3in][l]{{\footnotesize 3.}} 
                   {\it fold} the graph obtained so far;  

\smallskip
 
\makebox[.3in][l]{}
{\tt (*} The folded graph constructed so far, with all hairs ignored,
is called a  

\makebox[.3in][l]{}
``partial Cayley graph''. {\tt *)}  
% }    %parbox

\bigskip

\bigskip

\noindent We call {\it TC} a ``process'' (as opposed to ``algorithm'') 
because it does not terminate. 

\begin{defn} \ 
A {\bf Todd-Coxeter radius} of a finite presentation
$\langle A; R \rangle$ of a group $G$ is any function
$\rho_{_{\rm TC}}: \mathbb{N} \to \mathbb{N}$ with the following property:

After some number of steps the Todd-Coxeter process {\it TC} constructs 
a partial Cayley graph, called {\rm TC}$_n$ with radius 
$\leq \rho_{_{\rm TC}}(n)$ such that 

\smallskip

$(\forall w \in (A^{\pm 1})^*, |w| \leq n)$ 
{\bf [} $w =_{_{\langle A; R \rangle}} 1$ \ iff \ 
 {\sf red}$(w)$ labels a loop at the origin in {\rm TC}$_n$ {\bf ]}.
\end{defn}
When $n$ is even, the condition \ 
$(\forall w, |w| \leq n)[w =_{_{\langle A; R \rangle}} 1$ iff 
{\sf red}$(w)$ labels a loop at the origin] \ is equivalent to the following: 

{\it Within radius $n/2$ from the origin, the graph {\rm TC}$_n$ is 
identical to 
the ball of radius $n/2$ of the Cayley graph of the presentation}.

\smallskip

The concept of the Todd-Coxeter radius function appears indirectly in
\cite{FloydHoareLyndon} (in the case when it is linear). I learned about 
it from Stuart Margolis and John Meakin \cite{MargMeak}. 

\smallskip

If there exists a computable function which is an upper bound on 
$\rho_{_{\rm TC}}(\cdot)$
then the process {\it TC} can be used to decide the word problem of the
presentation $\langle A; R \rangle$.

We will view TC$_n$ as a DFA (deterministic finite automaton), by taking 
the origin $v_0$ as start and accept state, and by symmetrizing the edges 
(as we did for $\Lambda_{D(n)}$). 
We will also view TC$_n$ as a 2-complex (which agrees with the 
Cayley 2-complex within radius $n/2$ when $n$ is even). The context
will tell us which one of the three TC$_n$'s we are talking about.

\begin{pro} \ 
For any word $w \in (A^{\pm 1})^*$ of length $\leq n$ we have:  
 
\smallskip

$w =_{_{\langle A; R \rangle}} 1$ \  iff \ 
 the DFA {\rm TC}$_n$ accepts {\sf red}$(w)$.    
\end{pro} 
{\bf Proof.} \ This follows immediately from the fact that 
$w =_{_{\langle A; R \rangle}} 1$ iff {\sf red}$(w)$ 
labels a loop in the Cayley graph. Moreover, the partial Cayley graph 
TC$_n$ coincides with the Cayley graph within radius $n/2$.  
 \ \ \ $\Box$

\bigskip

Another way to build a 2-complex in order to solve the word problem for 
words of length $\leq n$ is as follows:
For each $w \in (A^{\pm 1})^*$ of length $\leq n$ we consider all the van 
Kampen diagrams of $w$ of minimum folded diameter. We create a new vertex 
$v_0$, and attach the origins of all these van Kampen diagrams to $v_0$; 
now we have a connected 2-complex. Next, we fold this 2-complex. 
We call this the {\bf folded van Kampen 2-complex} for words of length 
$\leq n$, and denote it by fK$_n$.  

\bigskip

A third way to build a 2-complex in order to solve the word problem for 
words of length $\leq n$ is as follows:
We create a vertex $v_0$ (an origin). For every relator $r \in R$ and 
every word $x \in (A^{\pm 1})^*$ of length $\leq \lambda(n)$ 
(for a certain function $\lambda(n)$ to be determined soon), 
we create a loop with origin $v_0$, 
labeled by $r^x$; this loop bounds one face. 
We denote this 2-complex by LC$_n$. Next, we fold 
the 2-complex LC$_n$, and denote the resulting 2-complex by fLC$_n$. 
Finally,  we choose $\lambda(n)$ large enough, but minimal, 
such that in fLC$_n$ we have: \  
For every word $w$ of length $\leq n$, $w =_{_{\langle A; R \rangle}} 1$ 
iff {\sf red}$(w)$ labels a closed path through the origin. 

The function $\lambda(.)$ above is called ``folded loop-complex 
function''.  One notes immediately that for the minimum function 
$\lambda(\cdot)$ of the presentation we have: 

\smallskip 

 \ \ \ \ \ \ \ LC$_n = \Lambda_{\lambda(n)}$,  

\smallskip 

\noindent where $\Lambda_j$ is the loop complex introduced following 
Definition \ref{isodiametric}. We call f$\Lambda_{\lambda(n)}$ the 
{\bf folded loop-complex} for words of length $\leq n$.  

\medskip

\begin{pro} \ 
The Todd-Coxeter 2-complex {\rm TC}$_n$ is equal to the folded loop 
2-complex {\rm f}$\Lambda_{\lambda(n)}$, 
and contains the folded van Kampen 2-complex {\rm fK}$_n$ as a 
subcomplex. 

The minimum Todd-Coxeter radius function $\rho_{_{\rm TC}}(.)$, 
the minimum folded isodiametric function $d(.)$, 
the minimum folded loop-complex function $\lambda(.)$, 
and the minimum isodiametric function $D(.)$,
are the same. 
\end{pro}
{\bf Proof.} \ 
(1) \ TC$_n$ can be ``pulled apart'' into loops with labels $r^x$, with 
$r \in R$, $|x| \leq \rho_{_{\rm TC}}(n)$ \ 
(one loop per face of the complex TC$_n$). 
More precisely, the process of {\bf pulling a complex apart into loops}
goes as follows: \\   
$\bullet$ For each face $f$ in the complex, choose a path $p_f$ of length 
$\leq \rho_{_{\rm TC}}(n)$ from the face to the origin of TC$_n$. 
Let $x$ be the label of $p_f$ and let $r \in R$ be the label of 
the contour of $f$. \\   
$\bullet$ Create a new origin for the loop complex to be constructed. \\
$\bullet$ Repeat the following, for each face $f$ of TC$_n$ until all faces 
have been removed from TC$_n$:  
 
- \ Create a new path with label $x$, attached at the new origin; 

 \ \ \ at the other end of this path, attach a face with contour label $r$ 

 \ \ \ (so, viewed from the new origin, this path-and-loop has label $r^x$).  

- \  Remove the face $f$ from its place in the TC complex.  

\smallskip 

\noindent The process of pulling TC$_n$ apart can be viewed as the inverse of
the folding process; it is reversible at each step. Therefore, if these loops 
are folded up again, we recover TC$_n$.  So we have: \ 

\smallskip

 \ \ \ \ \ \ \ \ \ \ \ $D(n) \leq \rho_{_{\rm TC}}(n)$, \ 

\smallskip

\noindent  where $D$ is the {\it minimum} isodiametric function.

\medskip

\noindent  (2) \ \ \ On the other hand, suppose we take {\it all} possible 
loops with label 
$r^x$, for every $r \in R$ and every reduced word $x$ of length 
$\leq \rho_{_{\rm TC}}(n)$, and attach these loops to TC$_n$ at the origin, 
and fold. We claim that the 2-complex obtained is again TC$_n$. 
Indeed, in each added loop the path labeled by $x$ has length 
$\leq \rho_{_{\rm TC}}(n)$.
By the minimality of the Todd-Coxeter radius function $\rho_{_{\rm TC}}(.)$, 
the process TC glues on all relators within radius 
$\rho_{_{\rm TC}}(n)$ anyway; 
hence, all $r^x$ (with $x$ reduced) occur already in TC$_n$.

So, TC$_n$ can be built by taking the folded loop-complex of radius 
$\rho_{_{\rm TC}}(n)$.  Since both $\rho_{_{\rm TC}}(.)$ and $\lambda(.)$ 
are minimal, we conclude that TC$_n$ = f$\Lambda_{\lambda(n)}$, and 
$\rho_{_{\rm TC}}(n) = \lambda(n)$.

\medskip

\noindent  (3) \ \ \ Since TC$_n$ can be obtained by folding loops with 
labels $r^x$ with $|x| \leq \rho_{_{\rm TC}}(n)$ 
(as seen at the beginning of the proof) we conclude that 

\smallskip

 \ \ \ \ \ \ \ \ \ \ \ $\rho_{_{\rm TC}}(n) \leq D(n)$. 

\smallskip

\noindent   Indeed, by Lemma \ref{isodiamNFA}, $D(n)$ is the radius of 
an unfolded loop-complex which can be used to decide the word problem for 
all words of length $\leq n$. Hence, the process {\it TC} will decide the
word problem after reaching radius $D(n)$.  

\medskip

\noindent (4) \ \ \ vK$_n$ can be pulled apart into loops with labels 
$r^x$, with $r \in R$, $|x| \leq d(n)$. 
This process of pulling vK$_n$ apart is reversible at each 
step; therefore, if these loops are folded up again, we recover vK$_n$.

Since f$\Lambda_{\lambda(n)}$ has minimum radius, we conclude that \ 
$\lambda(n) \leq d(n)$.

\medskip

\noindent (5) \ \ \ 
At the same time, each folded van Kampen diagram of a word of length
$\leq n$ can be obtained by folding loops with labels $r^x$. If in a folded 
van Kampen diagram with minimal diameter, more loops  are attached (at the 
origin) and folded in, this does not shrink the diameter (since the diameter
is already minimum); hence all minimum-diameter folded van Kampen diagrams 
of words of length $\leq n$, as well as the folded van Kampen complex 
fK$_n$ are subcomplexes of f$\Lambda_{\lambda(n)}$. 
Hence we also have $d(n) \leq \lambda(n)$ (since subcomplexes of 
f$\Lambda_{\lambda(n)}$ cannot have a larger radius than 
f$\Lambda_{\lambda(n)}$).

Hence, combining this with (4) we obtain, \ $\lambda(n) = d(n)$.

\medskip

\noindent (6) \ \ \ We saw in (1) that $D \leq \rho_{_{\rm TC}}$, 
we saw in (2) that $\rho_{_{\rm TC}} = \lambda$, and we saw in (3) that 
$\rho_{_{\rm TC}} \leq D$. Hence, $D = \rho_{_{\rm TC}} = \lambda$. 
We saw in (5) that $\lambda = d$. 
 \ \ \  $\Box$

\medskip

Since the four functions $\rho_{_{\rm TC}}(.), \lambda(.), d(.), D(.)$ 
are the same, we will use $d(.)$ to denote all of them.
The fact that $d(.)$ and $D(.)$ are the same appears implicitly in the 
literature (e.g., in \cite{Papa} the definition of the folded isodiametric 
function is used for the ``isodiametric function'', without mention that 
this is not the usual definition).

\begin{thm} \  
For any finite presentation, the minimum filling length function $F(.)$ 
and the minimum {\em folded} filling length function $f(.)$, are the same. 
\end{thm}
{\bf Proof.} \ Recall the characterization of the filling length in terms 
of a rewriting system (see the definition of filling length and the 
subsequent characterizations). The same rewriting characterization applies 
to the folded filling length, based on the folded van Kampen diagram fK$_n$. 
Both the minimum filling length and the minimum folded filling length  
for a word $w$ are equal to the length of 
the longest intermediate word derived in the rewrite process from $w$ to 
$1$. Hence the two functions $F$ and $f$ are equal.  
 \ \ \ $\Box$

\bigskip

\noindent {\bf Remark:} Minimum-area van Kampen diagrams may have much 
larger area than their folded version.
For example, consider a finite presentation $\langle A; \{ r \} \rangle$ 
where $r$ is cyclically reduced (relative to $FG(A)$) and $|r| > 0$, and
consider the word $w = r^n$ (for any $n > 1$). 
Then the van Kampen diagram of $w$, consisting of $n$ positive 
(counter-clockwise) loops labeled by $r$, attached at the origin, 
has area $\geq n$. But the folded van Kampen diagram has only one face. 

\begin{thm}
\label{ineqs} \  
Let $P_{\rm min}, f_{\rm min}, d_{\rm min}$ be the minimum isoperimetric 
function, respectively the minimum filling length function, respectively 
the minimum isodiametric function of a finite presentation 
$\langle A; R \rangle$. Let $p_{\rm min}$ be the 
minimum folded isoperimetric function.  These filling functions are related 
as follows (where \ $c >1$ is a constant that depends on the presentation; 
the constant may be different in different parts of the Theorem).

\medskip

\noindent {\bf (1)} \ \ \ \ \ \ \ \ 
$d_{\rm min}(n) \leq \frac{1}{2} \, f_{\rm min}(n) $
$\leq c \cdot (P_{\rm min}(n) + n)$, \ \ and \ \ 
$p_{\rm min}(n) \leq P_{\rm min}(n)$. 

\medskip

\noindent {\bf (2)} \ \ \ \ \ \ \ \  
             $P_{\rm min}(n) \leq  n \, c^{c^{d_{\rm min}(n)}}$
  \ \ \ \ \ \ \ \ (Cohen's double exponential inequality)

\medskip

\noindent {\bf (3)} \ \ \ \ \ \ \ \ 
              $P_{\rm min}(n) \leq c^{f_{\rm min}(n)}$ 
  \ \ \ \ \ \ \ \ (Gromov, Gersten)

\medskip
 
\noindent {\bf (4)} \ \ \ \ \ \ \ \ 
      $f_{\rm min}(n) \leq c^{d_{\rm min}(n) +n}$
  \ \ \ \ \ \ \ \ (Gersten, Riley)

\medskip

\noindent {\bf (5)} \ \ \ \ \    \ \ \ 
                     $p_{\rm min}(n) \leq c^{d_{\rm min}(n)}$.
\end{thm}
{\bf Proof} (or references): \ For the proofs of the first two inequalities 
in (1) see \cite{GerstenRiley1}.
The inequality $p \leq P$ is obvious. The double exponential 
inequality (2) is due to Daniel Cohen \cite{Cohen}; 
Steve Gersten \cite{GerstenDoubleEx} gave another proof,
and Papasoglu \cite{Papa} adapted Gersten's proof to more general 2-complexes.
We will give another proof of the double exponential inequality in a later 
section. 
(3) is due to Gromov and Gersten (\cite{Gromov}, pp.\ 100-101). 
(4) was proved by S.M.\ Gersten and T.\ Riley (Thm.\ 3 in 
\cite{GerstenRiley1}). 

(5) The folded isoperimetric function $p(\cdot)$ is bounded by the number 
of faces in a folded van Kampen diagram with minimum diameter $d(n)$. 
Since every vertex in a folded van Kampen diagram has degree 
$\leq 2 \, |A|$, it follows that a folded van Kampen diagram of diameter 
$d(n)$ has at most $a^{d(n)}$ edges (for some constant $a > 1$ depending on 
$|A|$). Every face has a boundary of length $\leq m$, where 
$m = {\rm  max}\{|r| : r \in R\}$ (so, $m$ depends on the presentation but
not on $n$). Therefore the folded van Kampen diagram
has $< (a^{d(n)})^m$ different boundary edge-cycles, hence it has 
$< c^{d(n)}$ faces (for a constant $c$).
 \ \ \ $\Box$

\bigskip

As a consequence of (4), (2) and (3) in the above Theorem we have the 
following break-up of the double exponential inequality (2) into two steps, 
when \ $d_{\rm min}(n) \geq a \, n$ \ (for a constant $a > 0$): 
 \ \ $P_{\rm min}(n) \ \leq \ c^{f_{\rm min}(n)}$
                      $ \ \leq \ c^{C^{d_{\rm min}(n)}}$

In \cite{GerstenRiley1} Gersten and Riley use
a slightly weaker form of the double exponential inequality, namely, \  
$P_{\rm min}(n) \leq  c^{c^{d_{\rm min}(n)+ n}}$.
By using (2) above, we can improve (4) and the break-up of (2):

\medskip

\noindent {\bf (4')} \ \ \ \ \   \ \ \ \ \ 
 $f_{\rm min}(n) \leq c^{d_{\rm min}(n)} + d_{\rm min}(n) \, \log n$  

\medskip

\noindent When \ $d_{\rm min}(n) \geq \log_a n$ \ (for a constant $a > 1$):
 
\smallskip

\noindent {\bf (2')} \ \ \ \ \   \ \ \ \ \
 $P_{\rm min}(n) \ \leq \ c^{f_{\rm min}(n)}$
   $  \ \leq \ c^{C^{d_{\rm min}(n)}}$.
 
\bigskip

\noindent {\bf Question:} How are the minimum filling length function
$f_{{\rm min}}$ and the minimum folded isoperimetric function
$p_{{\rm min}}$ related? 
Do we have \ $f_{{\rm min}}(n) \leq c \cdot (p_{{\rm min}}(n) + n)$ \ 
(for some constant $c > 0$)?

\bigskip

See \cite{GerstenRiley1}, \cite{GerstenRiley2}, and
\cite{Kapo}  for recent applications of isodiametric and other functions. 

Earlier we discussed the folded van Kampen diagrams, and we used them
to define the folded isoperimetric, isodiametric, and filling length
functions (the latter two were later shown to be equal to their unfolded
counterparts). We can define a further contraction of van Kampen diagrams
by mapping van Kampen diagrams into the Cayley 2-complex; let's call the
image of such a mapping of a van Kampen diagram the ``{\it Cayley image}
of the van Kampen diagram''. We can then define new functions:
The {\it Cayley isoperimetric function} (an upper bound on the number
of faces in the Cayley image of van Kampen diagrams for words of length
$n$), the {\it Cayley isodiametric function} (an upper bound on the
diameter of the Cayley image), and the {\it Cayley filling length function}
(an upper bound on the length of the homotopy loop within the Cayley complex,
as word of length $n$ is contracted to a point).

This kind of mapping of van Kampen diagrams is different than folding;
in the folding process we identify vertices (of the van Kampen diagram) that
are equivalent in the free group; in the Cayley map, we identify
vertices (of the van Kampen diagram) that are equivalent modulo the
group $G$ under consideration.

%%%%%%%%%%%%%%%%%%%%%%%%%%%%%%%%%%%%%%%%%%%%%%%%%%%%%%%%
% Section 3
%%%%%%%%%%%%%%%%%%%%%%%%%%%%%%%%%%%%%%%%%%%%%%%%%%%%%%%%
\section{Linear compression of the isoperimetric function}

Computational complexity is usually studied up to big-O because of the
linear speed-up theorem and the linear space compression theorem (see 
\cite{HU} for a reference). The filling functions are algebraic invariants 
up to big-O too; in addition, below we give an analogue of the linear 
speed-up and compression theorems for the isoperimetric function. 
It is not clear whether such a compression is possible for the 
isodiametric function and the filling length function. 

\begin{thm} \ 
Let $G$ be a group that has a finite presentation with isoperimetric 
function $\leq P(\cdot)$. Then $G$ also has a 
finite presentation with respect to which the isoperimetric function is \ 
$\leq P(n)/2 \ + \ n/2$. 
\end{thm}
{\bf Proof.} \ Let $\langle A; R \rangle$ be a finite presentation of $G$ 
with respect to which the isoperimetric function is $\leq P(\cdot)$. Let
$m =$ max$\{|r| : r \in R\}$ (length of the longest relator in $R$). 

A new presentation of $G$ is obtained as follows. First, we symmetrize $R$, 
i.e., for each $r \in R$, we add $r^{-1}$ and all cyclic permutations of $r$ 
and of $r^{-1}$ as relators. Let $\langle A; R_s \rangle$ be the symmetrized 
presentation obtained; this is still a presentation of the group $G$, with the 
same number $m$ defined above. 
Second, for any $i$-tuple of relators $(r_1, \ldots, r_i) \in (R_s)^i$ (with 
$2 \leq i \leq m$), we introduce the new relator {\sf red}$(r_1 \ldots r_i)$. 
In terms of van Kampen diagrams this means that we glue $r_1, \ldots, r_i$ 
together along a part of their boundaries, starting at the origins of the 
relators. We call the set of newly created relators $R_2$.
Obviously, $\langle A; R_s \cup R_2 \rangle$ is a finite presentation of $G$. 

\smallskip

We claim that the isoperimetric function of 
$\langle A; R_s \cup R_2 \rangle$ is $\leq P(n)/2 + n/2$. 

\smallskip

For a word $w \in (A^{\pm 1})^*$ with $n = |w|$, if $w =_G 1$ then there 
is a van Kampen diagram $K$ (over the original presentation 
$\langle A; R \rangle$) of area $\leq P(n)$. Let $K^*$ be the dual graph 
of $K$, and let $T^*$ be a spanning tree (a.k.a.\ maximal subtree) of $K^*$, 
whose root
is chosen to be the outer (unbounded) face. Let us now remove the root of 
$T^*$; this yields a forest $F^*$, with $\leq P(n)$ vertices. 
For each member tree of $F^*$ we choose the child of the root of $T^*$ as 
the root. 
Since the root of $T^*$ has degree $n$ in $K^*$, there are $\leq n$ member 
trees in the forest $F^*$. 

\smallskip

We will now use $F^*$ to transform the van Kampen diagram $K$ (over the 
presentation $\langle A; R \rangle$) into a van Kampen diagram of area \   
$\leq P(n)/2 + n/2$ \  over the new presentation 
$\langle A; R_s \cup R_2 \rangle$. The main observation is that each vertex
in $F^*$ has degree $\leq m$, and each tree root of $F^*$ has degree
$\leq m-1$ in $F^*$.

\smallskip

\noindent {\bf 1.} Let $n_1$ ($ \leq n$) be the the number of member trees 
in the forest $F^*$ that consist on only one vertex. 
We leave that part of $F^*$ alone.  

\smallskip

\noindent {\bf 2.} For each member tree of $F^*$ that has at least two 
vertices we do the following. We consider a maximal set $S$ of sibling 
leaves at maximum depth (siblings are vertices with the same parent). 
We fuse all the siblings in $S$ and their parent, thus forming a new vertex. 
In the van Kampen diagram $K$, this corresponds to fusing $\leq m$ 
neighboring faces into one new face over the new presentation 
$\langle A; R_s \cup R_2 \rangle$.
From now on we ignore this new vertex (remove it from the picture). 

\smallskip

\noindent {\bf 3.} We repeat step {\bf 2} as often as possible. When we reach 
the root of a member tree of $F^*$, either it still has children (which are
leaves now); then we fuse the root with these children into a new vertex. 
Or all the children were already removed.
In the latter case, we fuse the root with any one of the new vertices that
a child belongs to. Since every tree root in $F^*$ has degree $\leq m-1$, 
this creates a new vertex of $\leq m$ old vertices. 
 
\smallskip
 
As a result, we obtain a van Kampen diagram over 
$\langle A; R_s \cup R_2 \rangle$
with the following upper bound on the number of vertices:

\smallskip

\noindent 
\makebox[1.5in][l]{$\leq n_1$}  
             (for the $n_1$ one-vertex member trees of the forest $F^*$) \\ 
\makebox[1.5in][l]{$+ \ (P(n) - n_1)/2$} 
             (for the multi-vertex trees of the forest $F^*$, in which each 
              vertex was \\ 
\makebox[1.5in][l]{}  fused with at least one other vertex; and at least one 
              vertex was \\ 
\makebox[1.5in][l]{}  fused with more than one other vertex) \\  
$\leq n_1/2 \ + \ P(n)/2 \ \ \leq \ \ n/2 \ + \ P(n)/2$.
  \ \ \ \ \ $\Box$

%%%%%%%%%%%%%%%%%%%%%%%%%%%%%%%%%%%%%%%%%%%%%%%%%%%%%%%%
% Section 4 
%%%%%%%%%%%%%%%%%%%%%%%%%%%%%%%%%%%%%%%%%%%%%%%%%%%%%%%%
\section{A proof of the double exponential inequality, based on 
context-free languages }

The double exponential inequality gives an upper bound on the minimum 
isoperimetric function $P(\cdot)$ in terms of the minimum isodiametric 
function $d(\cdot)$. It is surprising that $d$ should provide any bound 
at all on $P$.  

There are many similarities between 
combinatorial group theory and the ``low-complexity'' theory of computation 
(see for example, \cite{AvenhausMadlenerI}, \cite{AvenhausMadlenerII}, 
\cite{AvenhausMadlener3},\cite{AvenhausMadlener4}, \cite{AvenhausMadlener5},
 \cite{Bi}, \cite{BORS}, \cite{CohenMadlenerOtto}, \cite{Epst},
\cite{MO}, \cite{SBR},  \cite{HUSimon}, \cite{Valiev}, \cite{Waack}). 
An interesting consequence of the following proof is 
that, from the point of view of the theory of computation, the double 
exponential inequality belongs into the theory of context-free languages.

\begin{thm} \ 
If $P(\cdot)$ is the minimum isoperimetric function of a finite presentation
$\langle A; R \rangle$ and $d(\cdot)$ is the minimum isodiametric function 
of that presentation then we have for all $n$: \ 
$$P(n) \leq  \ n \, 2^{C \, c^{d(n)}}$$
where $C = 2 \,(2|A| + 1) \, \|R\|^2$ and $c = (2|A|)^2$. 
\end{thm}
{\bf Proof.} \ We will use Lemma \ref{isoperimTreeNFA}.
We fix a word $w \in (A^{\pm 1})^*$, 
of length $|w| = n > 0$, and we assume $w =_{_{\langle A; R \rangle}} 1$.
We consider the reduced word \ {\sf red}$(w) = v_1 v_2 \ldots v_m$ \ with \ 
$|{\sf red}(w)| = m \leq n$.

\medskip

\noindent {\bf (1) \ } 
The language $L({\it tree}\Lambda_{d(n)})$ is of course a regular language,
accepted by the NFA {\it tree}$\Lambda_{d(n)}$, which has \ 
$\leq \|R\| \, (2|A|)^{d(n)}$ states, as we saw just before Lemma 
\ref{isoperimTreeNFA}. Notation: The set of next states of 
{\it tree}$\Lambda_{d(n)}$, reached from state $q$ under input letter
$a$, will be denoted by $q\cdot a$. The accept (and start) state of 
{\it tree}$\Lambda_{d(n)}$ will be denoted by ${\bf 1}_{\Lambda}$.
 
It is well known that $[w]_{{\rm{FG}}}$ is a context-free language, 
accepted by a push-down automaton (a ``pda'') with $O(n)$ states
(see \cite{HU}, \cite{Harrison} for background on 
{\it Dyck languages}, and on context-free languages in general). 

Here is a more detailed description of this pda, $\Pi_w$. The state set is 
$Q = \{f, s_0, s_1, \ldots, s_m\}$, where 
$s_m$ is the start state and $f$ is the accept state. 
Recall that $m = |{\sf red}(w)| \leq n$. The stack alphabet 
is $\Gamma = A^{\pm 1} \cup \{ {\bf z} \}$, where {\bf z} is the bottom 
marker of the stack and is also the initial content of the stack.
The input alphabet is $\Sigma = A^{\pm 1}$. 
As before, we will denote the empty word by 1 (in \cite{HU} it is denoted 
by $\varepsilon$). The transition relation \  
$\delta : Q \times (\Sigma \cup \{ 1 \}) \times \Gamma \ $
$\to \ Q \times \Gamma^*$ \ is defined as follows:

\smallskip

\noindent {\sf Phase 1}: Pop the top letter off the stack if the next 
input letter is the inverse of the top of the stack; otherwise, push 
the input letter on top of the stack.

$\delta(s_m, a, a^{-1}) = (s_m, 1)$, \ for all $a \in \Sigma$; 

$\delta(s_m, a, b) = (s_m, ab)$, \ for all $a \in \Sigma$,
 $b \in \Gamma$ \ with $b \neq a^{-1}$.

\smallskip

\noindent {\sf Phase 2}: Guess that the input is finished. Now, using 
``empty-input moves'', pop the stack and check that its content is the 
fixed word {\sf red}$(w)$
(with the beginning of the word at the bottom of the stack):

$\delta(s_i, 1, v_i) = (s_{i-1}, 1)$ \ for $i = m, \ldots, 1$;  

$\delta(s_0, 1, {\bf z}) = (f, 1)$.
\smallskip

\medskip

\noindent {\bf (2) \ } 
The intersection of a regular language and a context-free language is a 
context-free language; a pda for the intersection can be obtained thanks to 
a {\it cartesian product construction} (see \cite{HU} Theorem 6.5, or 
\cite{Harrison} Theorem 6.4.1). Let $\Pi$ be this pda accepting \ 
$[w]_{{\rm{FG}}} \cap L({\it tree}\Lambda_{d(n)})$,
obtained by the cartesian product construction.

Let us describe the pda $\Pi$ in more detail.
The stack alphabet of $\Pi$ is $\Gamma = A^{\pm 1} \cup \{ {\bf z} \}$, 
and the input alphabet is $\Sigma = A^{\pm 1}$, as before. For the 
state set of $\Pi$ we could take the cartesian product of the state set of 
${\it tree}\Lambda_{d(n)}$ and the state set of $\Pi_w$, 
but we can leave out the states that will not occur in any accepting 
computation. Hence, the states we keep form 
the set $Q = Q_1 \cup Q_2$ where \ $Q_1 = Q_{\Lambda} \times \{ s_m \}$  \  
and \  $Q_2 = \{{\bf 1}_{\Lambda}\} \times \{s_{m-1}, \ldots, s_1, s_0, f \}$. 

\smallskip
 
The transitions of $\Pi$ form two groups (as in the case of 
$\Pi_w$), which we call phase 1 and phase 2.
The first subset $Q_1$ corresponds to phase 1, and has
$|Q_1| = |Q_{\Lambda}| \leq \|R\| \cdot (2|A|)^{d(n)}$ states.
In phase 1 the transitions are 

\smallskip

$\delta((q, s_m), a, a^{-1}) = \{((p,s_m), 1) : p \in q\cdot a \}$, \ \ for 
      $a \in \Sigma, \ q \in Q_{\Lambda}$;

\smallskip

$\delta((q, s_m), a, b) = \{((p, s_m), ab) : p \in q\cdot a \}$, \ \ for 
   $b \in \Gamma, \ a \in \Sigma, \ q \in Q_{\Lambda}$, with $b \neq a^{-1}$.

\smallskip

\noindent The second subset $Q_2$ is used in phase 2, and has 
$m+1$ states. In phase 2 the transitions are 

\smallskip

$\delta(({\bf 1}_{\Lambda}, s_i), 1, v_i) = (({\bf 1}_{\Lambda}, s_{i-1}), 1)$,
 \ $i = m, \ldots, 1$;

\smallskip

$\delta(({\bf 1}_{\Lambda}, s_0), 1, {\bf z}) = (({\bf 1}_{\Lambda}, f), 1)$.

\smallskip

\noindent It is important to note that the $m+1$ states in $Q_2$ appear only 
in pop moves. 

The start state is $({\bf 1}_{\Lambda}, s_m)$ and the accept state is
$({\bf 1}_{\Lambda}, f)$. When the pda $\Pi$ reaches its accept
state its stack will always become empty; so, $\Pi$ ``accepts by empty stack''.

\medskip

\noindent {\bf (3) \ } 
Next, from the pda $\Pi$ (accepting by empty stack) we construct a 
context-free grammar that generates the language \  
$[w]_{{\rm{FG}}} \cap L({\it tree}\Lambda_{d(n)})$, 
thanks to a {\it construction of Chomsky, Evey, and Sch\"utzenberger} 
(see \cite{HU} Section 5.3, or \cite{Harrison} Theorem 5.4.3).  

The set of non-terminals corresponding to phase 1 are $S$ (the start
symbol of the grammar), and all symbols of the form 
$[p, c, q] \in Q \times \Gamma \times Q$.
The rules of phase 1 are of the form 

\smallskip

$S \to [({\bf 1}_{\Lambda}, s_m), {\bf z}, ({\bf 1}_{\Lambda}, f)]$;

\smallskip

$[(q, s_m), a, (p,s_m)] \to a^{-1}$, \ for any $(q, s_m) \in Q_1$, and 
   $p \in q\cdot a$;

\smallskip

$[(q,s_m), b, r_2] \to a [(p,s_m), a, r_1] [r_1, b, r_2]$, \ \

 \ \ \ \ \ \ \ \ \ for any $r_1, r_2 \in Q$ and  $(q,s_m), (p,s_m) \in Q_1$ 
   with \ $((p, s_m), ab) \in \delta((q, s_m), a, b)$.

\smallskip

\noindent The set of non-terminals corresponding to phase 2 is \   

$\{ [({\bf 1}_{\Lambda}, s_0), {\bf z}, ({\bf 1}_{\Lambda}, f)] \} \ \cup \ $
$\{ [({\bf 1}_{\Lambda}, s_i), v_i ,({\bf 1}_{\Lambda}, s_{i-1})] : $ 
$i = 1, \ldots, m \}$. \ 

\noindent These non-terminals belong to $Q_2 \times \Gamma \times Q_2$, 
except for $[({\bf 1}_{\Lambda}, s_m), v_m, ({\bf 1}_{\Lambda}, s_{m-1})]$, 
which belongs to $Q_1 \times \Gamma \times Q_2$. 
The only rules that have these non-terminals on the left-side are 
``empty-word rules''

\smallskip

$[({\bf 1}_{\Lambda}, s_i), v_i ,({\bf 1}_{\Lambda}, s_{i-1})] \ \to \ 1$, 
 \ for $i = 1, \ldots, m$; 

\smallskip

$[({\bf 1}_{\Lambda}, s_0), {\bf z}, ({\bf 1}_{\Lambda}, f)] \ \to \ 1$.

\medskip

\noindent {\bf (4) \ } We can simplify our grammar. First, we
drop the non-terminals of phase 2 altogether, since they only generate the 
empty word; we directly replace them by the empty word wherever they occur 
in the grammar. Thus, we assume from now on that our grammar contains no 
non-terminals in $Q_2 \times \Gamma \times Q_2$. 

We can also drop $S$ and use 
$[({\bf 1}_{\Lambda}, s_m), {\bf z}, ({\bf 1}_{\Lambda}, f)] $ 
$(\in Q_1 \times \Gamma \times Q_2)$ asthe start symbol.

\smallskip

We can discard all non-terminals in $Q_2 \times \Gamma \times Q_1$ 
because such non-terminals do not occur on the left side of any rule of the 
grammar. As a consequence, in every rule of the form 

$[(q,s_m), b, r_2] \to a [(p,s_m), a, r_1] [r_1, b, r_2]$

\noindent we now have \ $r_1 \in Q_1$. Hence, non-terminals in 
$Q_1 \times \Gamma \times Q_1$ generate only non-terminals that are also in
$Q_1 \times \Gamma \times Q_1$. On the other hand, the words generated, in one
step, by non-terminals in $Q_1 \times \Gamma \times Q_2$  are in \ 
$\Sigma \, (Q_1 \times \Gamma \times Q_1) \, (Q_1 \times \Gamma \times Q_2)$.
The main consequence of this is:

\smallskip

\noindent {\bf Fact:} \   
{\it In a parse tree of a word in $\Sigma^*$, only the right-most path can 
contain non-terminals in $Q_1 \times \Gamma \times Q_2$. 
All other non-terminals in the parse tree belong to 
$Q_1 \times \Gamma \times Q_1$.}
 
\medskip

\noindent In the following we will need bounds on the number of 
non-terminals (recall that \ 
$|Q_1| = |Q_{\Lambda}| \leq \|R\| \, (2|A|)^{d(n)}$, \    
$|Q_2| = m+1 \leq n+1$, and  $|\Gamma| = 2 |A| + 1$):  

\smallskip

$|Q_1 \times \Gamma \times Q_1| \ \leq \ (2|A| + 1) \|R\|^2 \, (2|A|)^{2d(n)}$
$ \ = \ C_1 \, c^{d(n)}$, 

\smallskip

$|Q_1 \times \Gamma \times Q_2| \ \leq \ (n+1) (2|A| + 1) \|R\| \, (2|A|)^{d(n)}$
$\ (<  \ n \, C_1 \, c^{d(n)})$, 

\smallskip

\noindent where \ $C_1 = (2|A| + 1) \, \|R\|^2$, and $c = (2|A|)^2$. 

\bigskip

\noindent {\bf (5) \ } 
We now use the {\it Pumping Lemma} (due to Bar-Hillel, Perles, Shamir,
see \cite{HU} Section 6.1, or \cite{Harrison} Theorem 6.2.1 and Corollary) 
which, among other things, states the following: 
If a language $L \neq \emptyset$ has a context-free grammar with $\nu$ 
non-terminals then $L$ contains a word of length $\leq \ell^{\nu}$, where 
$\ell$ is the maximum length of the right side of any rule. In our grammar,
$\ell = 3$ and $\nu \leq n c^{d(n)}$ for some constant $c > 1$. 
Thus, we immediately get an upper bound $3^{n c^{d(n)}}$ on the length of the 
shortest word. However, we can obtain a smaller upper bound if we use the above
Fact in our analysis of parse trees.

Recall that the Pumping Lemma is proved by looking at recurrences of non-terminals
on any path of the parse tree. A shortest word in the language will have a 
parse tree with no recurrent non-terminals on any path from the root. 
Hence, the right-most path of the parse tree has length \ 

\smallskip

$\leq |Q_1 \times \Gamma \times Q_2 \ \ \cup \ \ Q_1 \times \Gamma \times Q_1| $

\smallskip

$ \leq \ (n+1) (2|A| + 1) \|R\| \, (2|A|)^{d(n)} \ + \ $
$(2|A| + 1) \|R\|^2 \, (2|A|)^{2d(n)} $

\smallskip

 $\leq n \, C_1 \, c^{d(n)}$, \ \ \ with $c = (2 |A|)^2$ and 
        $C_1 = (2 |A| + 1) \, \|R\|^2$. 

\smallskip

\noindent By the Fact above, elsewhere the non-terminals in the parse tree 
that are not on the right-most path belong to 
$\in Q_1 \times \Gamma \times Q_2$. In other words, every subtree of the parse 
tree, hanging at a vertex of the right-most path, has only 
non-terminals in $Q_1 \times \Gamma \times Q_1$. 
Since there are no recurrent non-terminals, each one of these subtrees 
has depth \ 

\smallskip

$\leq |Q_1 \times \Gamma \times Q_1| \ \leq \ $
  $(2|A| + 1) \, \|R\|^2 \, (2|A|)^{2d(n)} = C_1 \, c^{d(n)}$.

\smallskip

\noindent Hence, since each non-terminal has at most 2 non-terminal successors,
the number of non-terminals in such a subtree is \ 

\smallskip

$\leq 2^{C_1 \, c^{d(n)}}$.  

\smallskip

\noindent Since the number of these subtrees is equal to the length of the
right-most path, we find that the total number of non-terminals in the parse 
tree is  

\smallskip

$\leq \ n \, C_1 \, c^{d(n)} \cdot 2^{C_1 \, c^{d(n)}}$
$ \ \leq  \ n \, 2^{2 \, C_1 \, c^{d(n)}}$  \ \    
(using the fact that \ $C_1 \, c^{d(n)} < 2^{C_1 \, c^{d(n)}}$).

\smallskip

\noindent Every non-terminal vertex has at most one terminal descendant 
in the parse tree, so the above upper bound also bounds the length of the 
shortest word in the language.  
Thus we have proved that for all words $w$ of length $|w| \leq n$:  
 
\smallskip 
 
$w =_{_{\langle A; R \rangle}} 1$ \ \ iff \ \     
$[w]_{{\rm{FG}}} \ \cap \ L({\it tree}\Lambda_{d(n)})$  \  
contains a word of length  \ $\leq \ n \, 2^{C \, c^{d(n)}}$, 

\smallskip

\noindent where $C = 2 \, C_1$.   
Now the Theorem follows from Lemma \ref{isoperimTreeNFA}. 
 \ \ \ $\Box$

%%%%%%%%%%%%%%%%%%%%%%%%%%%%%%%%%%%%%%%%%%%%%%%%%%%%%%%%
% Section5 
%%%%%%%%%%%%%%%%%%%%%%%%%%%%%%%%%%%%%%%%%%%%%%%%%%%%%%%%
\section{Computational complexity}

We know from \cite{Bi} and \cite{SBR} that if $P(\cdot)$ is an 
isoperimetric function for a
finite presentation $\langle A; R \rangle$, then the word problem of 
$\langle A; R \rangle$ is in NTime$(P)$.

\begin{pro} \  
Let $\langle A; R \rangle$ be a finite presentation with isodiametric function 
$\leq d(\cdot)$. Then the word problem of $\langle A; R \rangle$ is in 
{\rm DTime}$(c^{d(\cdot)})$, where $c > 1$ is a constant depending on the 
presentation.
\end{pro}
{\bf Proof:} \ Given a word $w$ of length $n$, we can decide whether 
$w =_{_{\langle A; R \rangle}} 1$
as follows. Construct the NFA $\Lambda_{d(n)}$; this can be done 
deterministically in time $\leq c^{d(n)}$ (for some constant $c >1$).
Next, we fold this automaton  as in Lemma \ref{foldedDFA};
this can be done deterministically in time $\leq C^{d(n)}$
(for some constant $C >1$).
Finally, check whether the folded DFA accepts {\sf red}$(w)$. \ \ \  $\Box$

\bigskip

Tim Riley \cite{Riley} observed that if a finite presentation has a filling 
length function $f$ then the word problem of that presentation has 
nondeterministic space complexity $O(f)$.

The following proposition strengthens this fact, by using {\it symmetric} 
Turing machines. Those are nondeterministic Turing machines 
whose transition relation is symmetric (i.e., the reverse of any transition 
of the machine is also a transition of that machine); see \cite{LewisPapad}
and \cite{Bi}. 
One can define space complexity in relation to such machines: 
SymSpace$(S)$ is the set of all languages accepted by symmetric
Turing machines with space $\leq S(\cdot)$. For time-complexity it is known 
that \ SymTime$(T) =$ NTime$(T)$ \ (proved by Lewis and 
Papadimitriou \cite{LewisPapad}). For space, DSpace$(S) \subseteq$ 
SymSpace$(S) \subseteq$ NSpace$(S)$; there are reasons to suspect that 
DSpace$(S) \neq$ SymSpace$(S) \neq$ NSpace$(S)$, but this remains an open 
problem.

\begin{pro} \ 
Let $\langle A; R \rangle$ be a finite presentation of a group with filling 
length function $\leq f(\cdot)$. Then the word problem of $G$ is in 
{\rm SymSpace}$(f(\cdot))$.
\end{pro}
{\bf Proof:} \  We use the rewriting system characterization of the filling
length function. Note that this rewrite system (described after the 
definition of isoperimetric functions) is symmetric. 
A symmetric Turing machine can simulate this rewrite system. 
Since the longest words that occur in the rewrite process 
have length $\leq f(n)$, the space needed by the Turing machine is also 
$\leq f(n)$.       \ \ \  $\Box$

\bigskip

The relations between filling functions on groups and the complexity of the 
word problem of groups are summarized below. Following the standard notation
for complexity classes, we introduce classes of finite presentations of groups,
based on their filling functions.

\begin{defn} \ 
Consider any function $h : \mathbb{N} \to \mathbb{N}$. We define \ 
{\bf Isoper}$(h)$ \ to be the set of all groups that have finite presentations 
$\langle A; R \rangle$ whose minimum isoperimetric function 
$P_{\langle A; R \rangle}$ satisfies 

\smallskip
 
$P_{\langle A; R \rangle}(n) \ \leq \ c_1 \ h(c_2 \, n)$, \ \ for all 
 $n \geq c_3$. 

\smallskip

\noindent Here, $c_1, c_2, c_3$ are positive constants, depending on
$\langle A; R \rangle$, but not on $n$. 
We say, ``the minimum isoperimetric function is $\leq h$ up to big-O''. 

In a similar way we define the sets of finite presentations 
{\bf Isodiam}$(h)$ for the isodiametric function, {\bf Filllen}$(h)$
for the filling length function, and {\bf FIsoper}$(h)$ for the folded
isoperimetric function.

By {\rm NTime}$(q)$ we denote all languages accepted by nondeterministic 
Turing machines with time complexity $O(q(O(n)))$. More precisely, for 
an accepted input of length $\leq n$ the Turing machine has at least one 
accepting computation whose time is \ $\leq c_1 \ q(c_2 \, n)$, \ for all
$n \geq c_3$; here, $c_1, c_2, c_3$ are positive constants, depending on
the Turing machine.
 
In a similar way we define {\rm DTime}$(q)$ and {\rm SymSpace}$(q)$.
\end{defn}

An inclusion between a class of groups and a class of languages (for
example, ${\rm Isoper}(q) \subseteq {\rm NTime}(q)$), is defined to mean  
that every group in ${\rm Isoper}(q)$ has its word problem in 
${\rm NTime}(q)$.

In this notation, the inequalities in Theorem \ref{ineqs} and the inclusions
in the above Propositions lead to the following (where $\subset$ and 
$\cap$ denote non-strict left-to-right or top-to-bottom inclusion). 

\begin{thm} \ 
For any function \ $q : \mathbb{N} \to \mathbb{N}$ \ with \    
$q(n) \geq \log \log n$ \ we have, 

\bigskip

\parbox{7in}{
${\rm Isoper}(q) \ \ \subset \ \ {\rm Filllen}(q) \ \ \ \subset \ \ \ $
$ {\rm Isodiam}(q) \ \ \ \subset \ \ {\rm Filllen}(2^{O(q)}) \ \ \ $
$ \subset \ \ \ {\rm Isoper}(2^{2^{O(q)}})$

\medskip

$ \ \ \ \ \cap \makebox[.9in]{} \cap \makebox[1.1in]{} \cap $
$\makebox[1.1in]{} \cap \makebox[1.2in]{} \cap $

\medskip

${\rm NTime}(q) \ \subset \ {\rm SymSpace}(q) \ \subset \ {\rm DTime}(2^{O(q)})$
$ \ \subset \ {\rm SymSpace}(2^{O(q)}) \ \subset \ {\rm NTime}(2^{2^{O(q)}})$
} %parbox

\bigskip

\noindent Moreover, \ \ \   
${\rm Isoper}(q) \subset {\rm FIsoper}(q)$, \ \  
and  \ \ 
${\rm Isodiam}(q) \subset {\rm FIsoper}(2^{O(q)})$.

\end{thm}
{\bf Proof.} \ To prove ${\rm Isoper}(q) \subset {\rm Filllen}(q)$, observe 
that if a group has a finite presentation $\langle A; R \rangle$ with 
isoperimetric function $\leq q$ then $\langle A; R \rangle$ has a filling 
length $\leq q$ too (since the minimum filling length is $\leq$ the minimum 
isoperimetric function up to big-O, by 
(1) of Theorem \ref{ineqs}). Hence, every presentation in ${\rm Isoper}(q)$ 
is also in ${\rm Filllen}(q)$. 

The other inclusions follow from Theorem \ref{ineqs} in a similar way.
 \ \ \ $\Box$

\bigskip

We do not know whether any of the inclusions in the above theorem are strict. 
For the complexity classes, this is a well known open problem.
Along the lines of \cite{Bi}, \cite{BORS} and \cite{SBR} one could make the 
following conjecture.

\medskip 

\noindent {\bf Conjecture.} A finitely generated group $G$ has its word 
problem in SymSpace$(S)$ iff $G$ is embeddable in a finitely presented group 
$H$ whose filling length function is $O(S)$.

%%%%%%%%%%%%%%%%%%%%%%%%%%%%%%%%%%%%%%%%%%%%%%%%%%%%%%%%

\bigskip

\bigskip

\noindent {\bf Acknowledgements.} I would like to thank John Meakin, 
Stuart Margolis, Ilya Kapovich, and especially Tim Riley
for enlightning discussions.

%%%%%%%%%%%%%%%%%%%%%%%%%%%%%%%%%%%%%%%%%%%%%%%%%%%%%%%%%%%%%%%%%%%%%%%%%%%%%%%

\bigskip 

\noindent
{\it Jean-Camille Birget  \\  
 Dept.\ of Computer Science  \\ 
 Rutgers University - Camden  \\ 
 Camden, NJ 08102, USA \\  
 birget@camden.rutgers.edu
}

\end{document}